\documentclass[a4paper,11pt]{article}

\usepackage{latexsym}
\usepackage{amsmath,amsfonts,amssymb,amsthm}
\newtheorem{thm}{Th\'eor\`eme}[section]

\newtheorem{proposition}[thm]{Proposition}

\newtheorem{remarque}[thm]{Remarque}
\newtheorem{definition}[thm]{D\'efinition}
\input epsf

\title{La suite de Thue-Morse et la cat\'egorie
$\hbox{Rec}$}
\author{Roland Bacher
}
\date{}
\begin{document}
\maketitle


R\'esum\'e:

Cette note introduit la cat\'egorie $\hbox{Rec}({\mathbf K})$
des matrices \`a r\'ecurrence sur un corps ${\mathbf K}$. Ceci permet de
calculer certains d\'eterminants reli\'es \`a la suite de Thue-Morse.

Abstract:
{\bf The Thue-Morse sequence and the category $\hbox{Rec}$. }
We define the category $\hbox{Rec}({\mathbf K})$ of recurrence
matrices over a field ${\mathbf K}$ and use it for calculating 
determinants of Hankel matrices related to the Thue-Morse sequence.



\section{Introduction}
\label{}

La fonction de (Prouhet-)Thue-Morse $\tau:{\mathbb N}\longrightarrow
  {\mathbb N}$ compte la somme des chiffres $\tau(\sum_{j=0}^l \epsilon_j2^j)=
\sum_{j=0}^l \epsilon_j$ d'un entier binaire $\sum_{j=0}^l
\epsilon_j2^j\in{\mathbb N}$ (avec $\epsilon_0,\dots,\epsilon_l\in
\{0,1\}$). Pour $n\geq 1$, soit $H(n)$ la
matrice de Hankel d'ordre $n$ avec coefficients complexes
$h_{s,t}=i^{\tau(s+t)}\in\{\pm 1,\pm i\},0\leq s,t< n,$ associ\'es
\`a la s\'erie g\'en\'eratrice
$\prod_{k=0}^\infty (1+ix^{2^k})=\sum_{n=0}^\infty i^{\tau(n)}x^n$.
Soit encore $f:\{1,2,\dots\}\longrightarrow\{\pm 1\}$ 
la fonction du pliage r\'egulier d\'efinie r\'ecursivement
par $f(2^n)=1,n\in{\mathbb N}$ et 
$f(2^n+a)=-f(2^n-a)$ pour tout $a$ tel que $1\leq a<2^n$ (voir \cite{AS}
pour plus d'informations sur Thue-Morse et le pliage r\'egulier,
voir \cite{APWW} pour des r\'esultats apparent\'es). 

\begin{thm} \label{det} On a l'\'egalit\'e
$\det(H(n+1))=\prod_{k=1}^n (1+i\ f(k))\in {\mathbb Z}[i]$ 
pour tout $n\in{\mathbb N}$ .
\end{thm}

La preuve, qui consiste \`a calculer la d\'ecomposition
$LU$ de $H(2^n)$,
fait intervenir une curieuse alg\`ebre li\'ee aux suites
automatiques, aux automates finis et aux groupes correspondants.

On pourrait en fait montrer le d\'eveloppement en 
fraction continue de type Jacobi
$$\prod_{k=0}^\infty
(1+ix^{2^k})=\frac{1}{1-u_0x-v_1x^2\frac{1}{
1-u_1x-v_2x^2\frac{1}{1-\dots}}}$$
avec $u_n=(-1)^n\ i,\ n\geq 0$ et o\`u la suite $v_1,v_2,\dots$ 
est d\'efinie par 
$v_1=1+i,v_2=1,v_3=-i,v_4=i,v_5=1,v_6=-i,v_7=1$
et pour $v_n,\ n=2^l+a\geq 8, 0\leq a<2^l,\ l\geq 3$ r\'ecursivement par
$v_n=i$ si $a\in \{0,2^{l-1}+1=\frac{n+2}{3}\}$, 
$v_n=1$ si $a\in \{1,2^{l-1}=\frac{n}{3}\}$ et 
$v_n=v_a$ autrement.

Des r\'esultats similaires existent \'egalement pour les
suites
$\beta_1,\beta_2,\beta_2,\dots$ et $\gamma_2,\gamma_3,\gamma_4,\dots$
\`a valeurs dans $\{0,\pm 1,\pm i,\pm 1\pm i\}$ d\'efinies par
$\sum_{n=0}^{\infty}\beta_nx^n=\frac{(1-x)}{i-1}
\prod_{k=0}^\infty (1+ix^{2^k})$ et
$\sum_{n=0}^{\infty}\gamma_nx^n=\frac{(1-x^2)}{i-1}\prod_{k=0}^\infty
(1+ix^{2^k})$.
Les d\'eterminants des matrices de Hankel associ\'ees ne prennent que
les valeurs $\pm 1,\pm i$.  

\section{Les cat\'egories ${\mathbf K}^{\mathcal M}$
et $\hbox{Rec}({\mathbf K})$}
Pour deux entiers naturels $p,q\in {\mathbb N}$ donn\'es,
${\mathcal M}_{p\times q}$ 
d\'esigne l'ensemble des paires de mots $(U,W)$ de m\^eme
longueur $l=l(U)=l(W)$ avec $U\in\{0,\dots p-1\}^l,\ W\in
\{0,\dots,q-1\}^l$. L'ensemble ${\mathcal M}_{p\times q}$ est donc
simplement le mono\"{\i}de libre engendr\'e par les 
$pq$ paires de mots $(s,t),\ 0\leq s<p,0\leq t<q$ de longueur $1$
pour la loi de composition $(U,W)(U',W')=(UU',WW')$.
Dor\'enavant, ${\mathcal M}_{p\times q}^l$ d\'esigne les paires 
de mots de longueur $l$ dans ${\mathcal M}_{p\times q}$.
Pour ${\mathbf K}$ un corps commutatif, fix\'e dans la suite,
une fonction $A\in {\mathbf K}^{
{\mathcal M}_{p\times q}}$ de ${
{\mathcal M}_{p\times q}}$ dans ${\mathbf K}$ d\'efinit
une suite de matrices de tailles $p^l\times q^l,\ l\in{\mathbb N}$, 
car on peut interpr\'eter les valeurs prises par $A$ sur l'ensemble 
fini ${\mathcal M}_{p\times q}^l$ comme 
coefficients d'une matrice de taille $p^l\times
q^l$ dont les indices parcourent ${\mathcal
  M}_{p\times q}^l$.

Ceci sugg\`ere de d\'efinir le {\it produit matriciel}
$A\cdot B\in{\mathbf K}^{{\mathcal M}_{p\times q}}$ de
$A\in{\mathbf K}^{{\mathcal M}_{p\times r}}$ et 
$B\in{\mathbf K}^{{\mathcal M}_{r\times q}}$ de la mani\`ere usuelle
en posant
$$(A\cdot B)[U,W]=\sum_{V\in \{0,\dots,r-1\}^l} A[U,V]B[V,W]$$
pour l'\'evaluation de $A\cdot B$ en 
$(U,W)\in{\mathcal M}_{p\times q}^l$.
 
Dans la suite, ${\mathbf K}^{\mathcal M}$ d\'esigne la cat\'egorie
dont chaque objet est un espace vectoriel ${\mathbf K}^{
{\mathcal M}_{q\times
    1}}$ et s'identifie \`a l'ensemble ${\mathbf K}^{{\mathcal M}_q}$
des fonctions sur le mono\"{\i}de libre ${\mathcal
  M}_q=\{0,\dots,q-1\}$ \`a $q$ g\'en\'erateurs.
Les morphismes de ${\mathbf K}^{{\mathcal M}_{q\times 1}}$
vers ${\mathbf K}^{{\mathcal M}_{p\times 1}}$ sont les \'el\'ements
de l'espace vectoriel ${\mathbf K}^{{\mathcal M}_{p\times
    q}}$. 

Un \'el\'ement 
$(S,T)\in{\mathcal M}_{p\times q}$ d\'etermine une application
lin\'eaire $\rho(S,T)\in
\hbox{End}({\mathbf K}^{{\mathcal M}_{p\times q}})$ 
en faisant correspondre \`a $A\in {\mathbf K}^{
{\mathcal M}_{p\times q}}$ la fonction $
\rho(S,T)A$ donn\'ee par
$(\rho(S,T)A)[U,W]=A[US,WT]$. Le petit calcul
$$\begin{array}{l}
\displaystyle
\rho(S,T)\big(\rho(S',T')A\big)[U,W]=\rho(S',T')A[US,WT]\\
\displaystyle \quad =A[USS',TT']=
\rho(SS',TT')A[U,W]\end{array}$$ montre qu'on obtient un morphisme de
mono\"{\i}des $\rho: {\mathcal M}_{p\times q}\longrightarrow
\hbox{End}({\mathbf K}^{{\mathcal M}_{p\times q}})$ d'image
le {\it mono\"{\i}de de d\'ecalage} $\rho({\mathcal M}_{p\times
  q})$. 
\begin{definition} Un sous-espace ${\mathcal A} \subset {\mathbf K}^{{
\mathcal M}_{p\times q}}$ est {\it r\'ecursivement clos} s'il est 
invariant par 
$\rho({\mathcal M}_{p\times q})$.
La cl\^oture r\'ecursive ${\overline A}^{rec}$ d'un \'el\'ement
$A\in{\mathbf K}^{{\mathcal M}_{p\times q}}$ est le plus petit
sous-espace r\'ecursivement clos contenant $A$. De mani\`ere
\'equivalente, ${\overline A}^{rec}$ est 
\'egalement le sous-espace engendr\'e par
l'orbite $\rho({\mathcal M}_{p\times q})A$ de $A$. La {\it 
complexit\'e} de $A$ est la cardinalit\'e $a=\hbox{dim}({\overline
A}^{
rec})\in{\mathbb N}\cup \{\infty\}$ d'une base de ${\overline
A}^{rec}$. Un \'el\'ement $A\in{\mathbf K}^{{\mathcal M}_{p\times q}}$
de complexit\'e finie $\hbox{dim}({\overline A}^{rec})<\infty$
est une {\it matrice \`a r\'ecurrence}. 
\end{definition}
On v\'erifie facilement
que l'ensemble 
$$\hbox{Rec}_{p\times q}({\mathbf K})=\{A\in{\mathbf K}^{{\mathcal
    M}_{p\times q}}\ \vert \ \hbox{dim}({\overline
  A}^{rec})<\infty\}$$
des matrices \`a r\'ecurrence est un espace vectoriel 
r\'ecursivement clos.

\begin{proposition} \label{produit}
Le produit $A\cdot B$ de deux matrices a r\'ecurrence
$A\in\hbox{Rec}_{p\times r}({\mathbf K})$, $B\in 
\hbox{Rec}_{r\times q}({\mathbf K})$ est une matrice \`a r\'ecurrence.
\end{proposition}
{\bf Preuve:} On choisit des g\'en\'erateurs $A_1,\dots,A_a$ 
et $B_1,\dots,B_b$ de
${\overline A}^{rec}$ et ${\overline B}^{rec}$. L'identit\'e
$$(A_i\cdot B_j)[Us,Wt]=\sum_{v=0}^{r-1} \big((\rho(s,v)A_i)\cdot
(\rho(v,t)B_j)\big)[U,W],\ $$
$(U,V)\in {\mathcal M}_{p\times q}$, 
$(s,t)\in {\mathcal M}_{p\times q}^1$, 
montre que l'espace vectoriel engendr\'e par les $ab$ produits 
$A_i\cdot B_j,\ 1\leq i\leq a,1\leq
j\leq b$, est r\'ecursivement clos.\hfill$\Box$
\begin{definition} 
La cat\'egorie $\hbox{Rec}({\mathbf K})$ des matrices \`a r\'ecurrence
est la sous-cat\'egorie de ${\mathbf K}^{\mathcal
    M}$ n'ayant que des fl\`eches dans $\hbox{Rec}_{p\times
    q}({\mathbf K})$.
Les objets de  $\hbox{Rec}({\mathbf K})$
sont les espaces $\hbox{Rec}_{q\times 1}({\mathbf K})$ 
des {\it vecteurs \`a r\'ecurrence}. \end{definition}
\begin{remarque} L'espace vectoriel $\hbox{Rec}_{p\times
  q}(
{\mathbf K})$ est un anneau pour le produit fonctionnel
$AB[U,W]=A[U,W]$ $B[U,W]$  
car le plongement ``diagonal'' de ${\mathbf K}^{{\mathcal
      M}_{p\times q}}$ dans ${\mathbf K}^{
{\mathcal M}_{pq\times pq}}$ pr\'eserve la complexit\'e.

L'espace
$\hbox{Rec}_{p\times  q}({\mathbf K})$ est aussi un anneau
(non-commutatif si $pq>1$) pour le produit de convolution 
$$(A*B)[U,W]=\sum_{(U,W)=(U_1,W_1)(U_2,W_2)}A[U_1,W_1]B[U_2,W_2]$$
obtenu en identifiant ${\mathbf K}^{{\mathcal M}_{p\times q}}\sim
{\mathbf K}^{{\mathcal M}_{(pq)}}$ avec l'anneau des s\'eries formelles
en $pq$ variables non-commutatives (ceci r\'esulte de l'identit\'e
$\rho(s,t)(A*B)=(\rho(s,t)A)(B[\emptyset,\emptyset])+A*
(\rho(s,t)B)$ pour
$(s,t)\in{\mathcal M}_{p\times q}^1$).
\end{remarque}

Un sous-espace ${\mathcal A}\subset
\hbox{Rec}_{p\times q}({\mathbf K})$ r\'ecursivement clos 
est compl\`etement caract\'eris\'e par l'action de 
$\rho({\mathcal M}_{p\times q})$ sur ${\mathcal A}$ et par
les \'evaluations 
$A\longmapsto A[\emptyset,\emptyset]$ en $(\emptyset,\emptyset)
\in{\mathcal M}_{p\times q}^0$ pour $A\in{\mathcal A}$. Une matrice
\`a r\'ecurrence $A$ de complexit\'e $a$ peut donc se d\'ecrire
\`a l'aide de $a$ {\it valeurs initiales}
$(A_1[\emptyset,\emptyset],\dots,A_a[\emptyset,
\emptyset])\in
{\mathbf K}^a$ (pour $A_1=A,\dots,A_a\in\hbox{Rec}_{p\times q}({\mathbf
  K})$ une base de ${\overline A}^{rec}$)
et de $pq$ {\it matrices de d\'ecalage} (abusivement not\'ees)
$\rho(s,t)\in \hbox{End}(\oplus_{h=1}^a
{\mathbf K}A_h)$ 
d\'efinies par $\rho(s,t)A_j=\sum_{k=1}^a\rho(s,t)_{k,j}A_k$
et d\'ecrivant l'action du mono\"{\i}de de d\'ecalage par rapport 
\`a la base $A_1,\dots, A_a$ de $\overline{A}^{rec}$. 
Par dualit\'e, une telle {\it pr\'esentation minimale} permet de
calculer une \'evaluation $A[U,W]$ de 
$A=A_1\in\hbox{Rec}_{p\times q}({\mathbf K})$ en utilisant la formule
$$\left(\begin{array}{c}A_1[s_1\dots s_n,t_1\dots t_n]\\\vdots\\
A_a[s_1\dots s_n,t_1\dots
t_n]\end{array}\right)=\rho(s_n,t_n)^t\cdots \rho(s_1,t_1)^t
\left(\begin{array}{c}A_1[\emptyset,\emptyset]\\\vdots\\
A_a[\emptyset,\emptyset]\end{array}\right).$$

Soit $A[{\mathcal M}_{p\times q}^{\leq n}]\in{\mathbf
  K}^{{\mathcal M}_{p\times q}^{\leq n}}$ la restriction
de $A\in{\mathbf K}^{{\mathcal M}_{p\times q}}$ \`a l'ensemble fini
${\mathcal M}_{p\times q}^{\leq n}$ des mots de longueur au plus $n$
dans ${\mathcal M}_{p\times q}$. Pour ${\mathcal A}\subset 
{\mathbf K}^{{\mathcal M}_{p\times q}}$ un espace vectoriel,
la notation ${\mathcal A}[
{\mathcal M}_{p\times q}^{\leq n}]\subset{\mathbf K}^{{\mathcal M}_{p\times
    q}}$ d\'esigne le sous-espace vectoriel
\'evident obtenu par la projection $A\longmapsto A[{\mathcal M}_{
p\times q}^{\leq n}]$.
\begin{definition}
Le {\it niveau de saturation} d'un espace vectoriel ${\mathcal A}$
est le plus petit \'el\'ement $N\in{\mathbb N}\cup
\infty$ tel que la projection naturelle
${\mathcal A}[{\mathcal M}_{p\times q}^{\leq N+1}]\longrightarrow
{\mathcal A}[{\mathcal M}_{p\times q}^{\leq N}]$
est un isomorphisme.
\end{definition}
\begin{proposition} \label{propfond}
Soit ${\mathcal A}\subset \hbox{Rec}_{p\times q}({\mathbf
    K})$ un espace vectoriel r\'ecursivement clos de niveau de
  saturation fini $N<\infty$. Alors 
${\mathcal A}$ et ${\mathcal A}[{\mathcal M}_{p\times
  q}^{\leq N}]$ sont isomorphes.
\end{proposition}

{\bf Id\'ee de la preuve} Notant $K_l\subset{\mathcal A}$ le
noyau de la projection \'evidente ${\mathcal A}\longrightarrow
{\mathcal A}[{\mathcal M}_{p\times q}^{\leq l}]$, on a l'\'egalit\'e
$K_N=K_{N+1}$ qui implique $K_n=K_N=\{0\}$ pour tout $n\geq N$.\hfill
$\Box$

La proposition \ref{propfond} permet de construire des pr\'esentations
minimales de $A+B$ et $A\cdot B$ pour $A,B$ des matrices \`a 
r\'ecurrence convenables (donn\'ees par des pr\'esentations minimales)
en utilisant un nombre fini d'op\'erations dans des espaces vectoriels
de dimensions finies. Plus pr\'ecis\'ement, \'etant donn\'ees des
bases $A_1=A,\dots,A_a$ et $B_1=B,\dots,B_b$ de $\overline{A}^{rec}$
et $\overline{B}^{rec}$, le calcul du niveau de saturation de
${\mathcal C}=\sum {\mathbf K} A_i+\sum {\mathbf K}B_j$ permet
de d\'eterminer le sous-espace $\overline{A_1+B_1}^{rec}\subset 
{\mathcal C}$ et d'en donner une base. Pour le produit, on 
proc\`ede similairement avec $\overline{A_1\cdot B_1}^{rec}
\subset \sum {\mathbf K} A_i\cdot B_j$.

\begin{remarque} L'alg\`ebre $\hbox{Rec}_{p\times p}({\mathbf K})$ 
contient des \'el\'ements inversibles
(pour le produit matriciel) dans
${\mathbf K}^{{\mathcal    M}_{p\times p}}$ mais 
sans inverse dans $\hbox{Rec}_{p\times p}({\mathbf K})$.
Un exemple est la matrice \`a r\'ecurrence diagonale
d\'efinie par 
$A[U,W]=1+n$ si $U=W=s_1\dots s_n,\ n\in{\mathbb N}$, et 
$A[U,W]=0$ sinon.
\end{remarque}

\section{Id\'ee de la preuve du th\'eor\`eme \ref{det}}

On exhibe des \'el\'ements $L,U\in \hbox{Rec}_{2\times 2}$ (o\`u 
$L[{\mathcal M}_{2\times 2}^l], U[{\mathcal M}_{2\times 2}^l]$
sont respectivement une matrice triangulaire unipotente inf\'erieure
et une matrice triangulaire sup\'erieure) tels que
le produit $H=L\cdot U\in \hbox{Rec}_{2\times 2}$
est donn\'e par $H[s_1\dots s_n,t_1\dots t_n]=\prod_{j=1}^n
i^{s_j+t_j}$
pour $s_1\dots s_n,t_1\dots t_n\in \{0,1\}^n$. Une
inspection des ``coefficients diagonaux'' de $U$ termine alors la preuve.
(Pour le d\'eveloppement en fraction
continue de Jacobi, on proc\`ede similairement en calculant la 
matrice de Stieltjes associ\'ee).

Plus pr\'ecis\'ement, on montre que
la matrice $H$ ci-dessus admet la pr\'esen\-tation minimale
$H_1=H,H_2$ avec valeurs initiales $H_1[\emptyset,\emptyset]=1,
H_2[\emptyset,\emptyset]=i$ et matrices de d\'ecalage 
$$
\rho(0,0)=\left(\begin{array}{cc}1&i\\0&0\end{array}\right),
\quad
\rho(0,1)=\rho(1,0)=\left(\begin{array}{cc}0&-i\\1&1+i\end{array}\right),\quad
\rho(1,1)=\left(\begin{array}{cc}i&i\\0&0\end{array}\right).$$

Similairement, $L$ peut se d\'ecrire par rapport \`a la
base $L_1=L,L_2,L_3,L_4$ par la pr\'esentation minimale avec
valeurs initiales $(L_1,\dots,L_4)[\emptyset,\emptyset]=(1,i,1,0)$ et 
les matrices de d\'ecalage
$$
\rho(0,0)=\left(\begin{array}{cccc}
1&i&1&0\\0&0&0&0\\0&0&0&0\\0&0&0&0\end{array}\right),\quad
\rho(0,1)=\left(\begin{array}{cccc}
0&0&0&0\\0&0&0&0\\0&0&0&0\\0&1&0&1\end{array}\right),$$
$$
\rho(1,0)=\left(\begin{array}{cccc}
0&-i&-1+i&-i\\1&1+i&-i&1\\0&0&0&0\\0&0&0&0\end{array}\right),
\quad
\rho(1,1)=\left(\begin{array}{cccc}
0&0&0&0\\0&0&0&0\\1&1+i&1&i\\0&i&0&i\end{array}\right).$$
La matrice \`a r\'ecurrence $U$ est le produit 
matriciel $U=D\cdot L^t$ avec $L^t[V,W]=L[W,V]$ et $D\in\hbox{Rec}_{
2\times 2}({\mathbb C})$ diagonal. Une pr\'esentation minimale de
$D$ est donn\'ee par $D_1=D,D_2,D_3$,
$(D_1,D_2,D_3)[\emptyset,\emptyset]=(1,1+i,1+i)$ et les matrices
de d\'ecalage 
$$\begin{array}{c}
\displaystyle \rho(0,0)=\left(\begin{array}{ccc}1&0&0\\0&0&0\\0&1&1
\end{array}\right),\quad
\rho(0,1)=\rho(1,0)=\left(\begin{array}{ccc}0&0&0\\0&0&0\\0&0&0
\end{array}\right),\\
\displaystyle 
\rho(1,1)=\left(\begin{array}{ccc}0&2&0\\1&1&1\\0&-2&0
\end{array}\right).\end{array}$$
Une inspection des coefficients de $D$ termine la preuve.

\begin{remarque} 
Le r\'esultat du th\'eor\`eme \ref{det} 
semble \'egalement vrai pour les s\'eries g\'en\'eratrices
$\prod_{k=0}^{\infty}(1+\sigma_k\ ix^{2^k})$ avec
$\sigma_0,\sigma_1,\dots \in \{\pm 1\}$ arbitraires
en rempla\c cant la suite du pliage r\'egulier par la suite d'un pliage
g\'en\'eralis\'e d\'efinie par $f(2^k)=\sigma_{k+1}$ et
$f(2^k+a)=-f(2^k-a),1\leq a<2^k$.
%
\end{remarque}


\begin{thebibliography}{00}

\bibitem{AS} J.-P. Allouche, J. Shallit, Automatic sequences. 
Theory, applications, generalizations, Cambridge
University Press (2003).

\bibitem{APWW} J.-P. Allouche, J. Peyri\`ere, Z.-X. Wen, Z.-Y. Wen,
{\it Hankel determinants of the Thue-Morse sequence.}
Ann. Inst. Fourier {\bf 48}, No.1, 1-27 (1998).


\end{thebibliography}
\end{document}